\documentclass[12pt]{amsart}
\usepackage{amsfonts}
\usepackage{amsmath,amssymb,amsthm}

\newcommand{\R}{\mathbb R}

\newcommand{\E}{\mathbb E}
\newcommand{\C}{\mathbb C}

\renewcommand{\span}{\mathrm{span}}

\newtheorem{thm}{Theorem}[section]
\newtheorem{cor}[thm]{Corollary}

\newtheorem{prop}[thm]{Proposition}
\theoremstyle{definition}

\theoremstyle{remark}

\newcommand{\ds}{\displaystyle}

\begin{document}

\title[GEOMETRIC INTERPRETATION OF THE INVARIANTS OF A SURFACE IN $\R^4$]
{GEOMETRIC INTERPRETATION OF THE INVARIANTS OF A SURFACE IN $\R^4$ VIA THE
TANGENT INDICATRIX AND THE NORMAL CURVATURE ELLIPSE}

\author{Georgi Ganchev and Velichka Milousheva}
\address{Bulgarian Academy of Sciences, Institute of Mathematics and Informatics,
Acad. G. Bonchev Str. bl. 8, 1113 Sofia, Bulgaria}
\email{ganchev@math.bas.bg}
\address{Bulgarian Academy of Sciences, Institute of Mathematics and Informatics,
Acad. G. Bonchev Str. bl. 8, 1113, Sofia, Bulgaria}
\email{vmil@math.bas.bg}



\begin{abstract}
At any point of a surface in the four-dimensional Euclidean space we consider
the geometric configuration consisting of two figures: the tangent indicatrix,
which is a conic in the tangent plane, and the normal curvature ellipse.
We show that the basic geometric classes of surfaces in the four-dimensional
Euclidean space, determined by conditions on their invariants, can be interpreted
in terms of the properties of the two geometric figures. We give some non-trivial
examples of surfaces from the classes in consideration.
\end{abstract}

\maketitle

\section{Introduction}

In this paper we deal with the theory of surfaces in the four-dimensional Euclidean
space $\R^4$.

Let $M^2$ be a surface in $\R^4$ with tangent space $T_pM^2$ at any point $p \in M^2$.
In \cite{GM1} we introduced the linear map $\gamma$ of Weingarten type at any $T_pM^2$
and sketched out the invariant theory of surfaces on the base of $\gamma$.

We show that the role of the map $\gamma$ in the theory of surfaces in $\R^4$ is similar
to that of the Weingarten map in the theory of surfaces in $\R^3$.

First, the map $\gamma$ generates two invariant functions $k$ and $\varkappa$,
analogous to the Gauss curvature and the mean curvature in $\R^3$. Here again
the sign of the function $k$ is a geometric invariant and the sign of $\varkappa$
is invariant under the motions in $\R^4$. However, the sign of $\varkappa$ changes
under symmetries with respect to a hyperplane in $\R^4$. The invariants $k$ and $\varkappa$
divide the points of $M^2$ into four types: flat, elliptic, hyperbolic and parabolic points.
In \cite {GM1} we gave a constructive classification of the surfaces consisting of flat
points, i.e. satisfying the condition $k= \varkappa = 0$. Everywhere, in the present
considerations we exclude the points at which $k= \varkappa = 0$.

Further, the map $\gamma$ generates the second fundamental form $II$ at any point
$p \in M^2$. The notions of a normal curvature of a tangent, conjugate and asymptotic
tangents are introduced in the standard way by means of $II$. The asymptotic tangents
are characterized by zero normal curvature.

The first fundamental form $I$ and the second fundamental form $II$ generate principal
tangents and principal lines, as in $\R^3$. Here, the points at which any tangent is
principal (''umbilical" points) are characterized by zero mean curvature vector, i.e.
the surfaces consisting of "umbilical" points are exactly the minimal surfaces in $\R^4$.
The principal normal curvatures $\nu'$ and $\nu''$ arise in the standard way and the
invariants $k$ and $\varkappa$ satisfy the equalities
$$k = \nu' \nu''; \qquad \varkappa = \frac{\nu' + \nu''}{2}.$$

The indicatrix of Dupin at an arbitrary (non-flat) point of a surface in $\R^3$ is
introduced by means of the second fundamental form. Here, using the second
fundamental form $II$, we introduce the indicatrix $\chi$ at any point $p \in M^2$
in the same way:
$$\chi : \nu' X^2 + \nu'' Y^2 = \varepsilon, \qquad \varepsilon = \pm1.$$
Then the elliptic, hyperbolic and parabolic points of a surface $M^2$ are characterized
in terms of the indicatrix $\chi$ as in $\R^3$. The conjugacy in terms
of the second fundamental form coincides with the conjugacy with respect to the
indicatrix $\chi$.

In \cite{GM1, GM2} we proved that the surface $M^2$ under consideration is with flat normal
connection if and only if $\varkappa =0$. In Section 3 we prove that:
\vskip 1mm
\emph{The surface $M^2$ is minimal if and only if the indicatrix $\chi$ is a circle.}
\vskip 1mm
\emph{The surface $M^2$ is with flat normal connection if and only if the indicatrix
$\chi$ is a rectangular hyperbola (a Lorentz circle).}

We also characterize the surfaces with flat normal connection in terms of
the properties of the normal curvature ellipse.

In Section 4 we give examples of surfaces with $\varkappa = 0$.

in Section 5 we give examples of surfaces with $k=0$.

\section{An interpretation of the second fundamental form}

Let $M^2: z = z(u,v), \, \, (u,v) \in {\mathcal D}$ (${\mathcal D}
\subset \R^2$) be a 2-dimensional surface in $\R^4$. The tangent
space $T_pM^2$ to $M^2$ at an arbitrary point $p=z(u,v)$ of $M^2$
is ${\rm span} \{z_u, z_v\}$.   We choose an orthonormal normal
frame field $\{e_1, e_2\}$ of $M^2$ so that the quadruple $\{z_u,
z_v, e_1, e_2\}$ is positive oriented in $\R^4$. Then the
following derivative formulas hold:
$$\begin{array}{l}
\vspace{2mm} \nabla'_{z_u}z_u=z_{uu} = \Gamma_{11}^1 \, z_u +
\Gamma_{11}^2 \, z_v
+ c_{11}^1\, e_1 + c_{11}^2\, e_2,\\
\vspace{2mm} \nabla'_{z_u}z_v=z_{uv} = \Gamma_{12}^1 \, z_u +
\Gamma_{12}^2 \, z_v
+ c_{12}^1\, e_1 + c_{12}^2\, e_2,\\
\vspace{2mm} \nabla'_{z_v}z_v=z_{vv} = \Gamma_{22}^1 \, z_u +
\Gamma_{22}^2 \, z_v
+ c_{22}^1\, e_1 + c_{22}^2\, e_2,\\
\end{array}$$
where $\Gamma_{ij}^k$ are the Christoffel's symbols and
$c_{ij}^k$, $i, j, k = 1,2$ are functions on $M^2$.

We use the standard denotations \;$E(u,v)=g(z_u,z_u), \;
F(u,v)=g(z_u,z_v), \; G(u,v)=g(z_v,z_v)$ for the coefficients of
the first fundamental form and set $W=\sqrt{EG-F^2}$. Denoting by
$\sigma$ the second fundamental tensor of $M^2$, we have
$$\begin{array}{l}
\sigma(z_u,z_u)=c_{11}^1\, e_1 + c_{11}^2\, e_2,\\
[2mm]
\sigma(z_u,z_v)=c_{12}^1\, e_1 + c_{12}^2\, e_2,\\
[2mm] \sigma(z_v,z_v)=c_{22}^1\, e_1 + c_{22}^2\,
e_2.\end{array}$$ \vskip 1mm

In \cite{GM1} we introduced a geometrically determined linear map
$\gamma$ in the tangent space at any point of a surface $M^2$ and
found invariants generated by this map.

We consider the functions
$$L = \frac{2}{W}\left|%
\begin{array}{cc}
\vspace{2mm}
  c_{11}^1 & c_{12}^1 \\
  c_{11}^2 & c_{12}^2 \\
\end{array}%
\right|, \quad
M = \frac{1}{W}\left|%
\begin{array}{cc}
\vspace{2mm}
  c_{11}^1 & c_{22}^1 \\
  c_{11}^2 & c_{22}^2 \\
\end{array}%
\right|, \quad
N = \frac{2}{W}\left|%
\begin{array}{cc}
\vspace{2mm}
  c_{12}^1 & c_{22}^1 \\
  c_{12}^2 & c_{22}^2 \\
\end{array}%
\right|.$$
Denoting
$$\displaystyle{\gamma_1^1=\frac{FM-GL}{EG-F^2}, \quad
\gamma_1^2 =\frac{FL-EM}{EG-F^2}}, \quad
\displaystyle{\gamma_2^1=\frac{FN-GM}{EG-F^2}, \quad
\gamma_2^2=\frac{FM-EN}{EG-F^2}},$$
we obtain the linear map
$$\gamma: T_pM^2 \rightarrow T_pM^2,$$
determined by the equalities
$$\begin{array}{l}
\vspace{2mm}
\gamma(z_u)=\gamma_1^1z_u+\gamma_1^2z_v,\\
\vspace{2mm} \gamma(z_v)=\gamma_2^1z_u+\gamma_2^2z_v.
\end{array}$$

The linear map $\gamma$ of Weingarten type at the point $p \in M^2$ is invariant
with respect to changes of parameters on $M^2$ as well as to motions in $\R^4$.
This implies that the functions
$$k = \frac{LN - M^2}{EG - F^2}, \qquad
\varkappa =\frac{EN+GL-2FM}{2(EG-F^2)}$$ are invariants of the surface $M^2$.

The invariant $\varkappa$ is the curvature of the normal connection of the surface $M^2$ in $\E^4$.

The invariants $k$ and $\varkappa$ divide the points of $M^2$ into four types \cite{GM1}:
flat, elliptic, parabolic and hyperbolic. The surfaces consisting of flat points
satisfy the conditions
$$k(u,v)=0, \quad \varkappa(u,v)=0, \qquad (u,v) \in \mathcal D,$$
or equivalently $L(u,v)=0, \; M(u,v)=0, \; N(u,v)=0, \; (u,v) \in \mathcal D.$
These surfaces are either planar surfaces (there exists a hyperplane $\R^3 \subset \R^4$
containing $M^2$) or developable ruled surfaces.

Further we consider surfaces free of flat points, i.e. $(L, M, N) \neq (0, 0, 0)$.
\vskip 2mm
Let $X = \alpha z_u + \beta z_v, \,\,(\alpha,\beta) \neq (0,0)$ be a tangent vector
at a point $p \in M^2$. The Weingarten map $\gamma$ determines a second fundamental
form of the surface $M^2$ at $p \in M^2$ as follows:
$$II(\alpha,\beta) = - g(\gamma(X),X) = L\alpha^2 + 2M\alpha\beta + N\beta^2, \quad \alpha,\beta \in \R.$$

As in the classical differential geometry of surfaces in $\R^3$ the second fundamental
form $II$ determines conjugate tangents at a point $p$ of $M^2$.

Two tangents $g_1: X = \alpha_1 z_u + \beta_1 z_v$ and $g_2: X = \alpha_2 z_u + \beta_2 z_v$
are said to be \textit{conjugate tangents}  if $II(\alpha_1, \beta_1; \alpha_2, \beta_2) = 0$,
i.e.
$$L\alpha_1 \alpha_2 + M (\alpha_1 \beta_2 +\alpha_2 \beta_1) + N\beta_1 \beta_2 = 0.$$

A tangent $g: X = \alpha z_u + \beta z_v$ is said to be \textit{asymptotic} if it is
self-conjugate, i.e.
$L\alpha^2 + 2M\alpha\beta + N\beta^2 = 0$.

A tangent $g: X = \alpha z_u + \beta z_v$ is said to be \textit{principal} if it is
perpendicular to its conjugate. The equation for the principal tangents at a point $p \in M^2$ is
$$\left|\begin{array}{cc}
E & F\\
[2mm] L & M \end{array}\right| \alpha^2+ \left|\begin{array}{cc}
E & G\\
[2mm] L & N \end{array}\right| \alpha \beta+
\left|\begin{array}{cc}
F & G\\
[2mm] M & N \end{array}\right| \beta^2=0.$$

A line $c: u=u(q), \; v=v(q); \; q\in J$ on $M^2$ is said to be a
\textit{principal line} (a \textit{line of curvature}) if its tangent at any point
is principal. The surface $M^2$ is parameterized with respect to the principal
lines if and only if
$$F=0, \qquad M=0.$$

Let $M^2$ be parameterized with respect to the principal lines and
denote the unit vector fields $\displaystyle{x=\frac{z_u}{\sqrt
E}, \; y=\frac{z_v}{\sqrt G}}$. The equality $M = 0$ implies that
the normal vector fields $\sigma(x,x)$ and $\sigma(y,y)$ are
collinear. We denote by $b$ a unit normal vector field collinear
with $\sigma(x,x)$ and $\sigma(y,y)$, and by $l$ the unit normal
vector field such that $\{x,y,b,l\}$ is a positive oriented
orthonormal frame field of $M^2$ (the two vectors $\{b, l\}$ are determined up to
a sign). Thus we obtain a geometrically determined orthonormal frame field
$\{x,y,b,l\}$ at each point $p \in M^2$. With respect to the frame field $\{x,y,b,l\}$
we have the following formulas:
$$\begin{array}{l}
\vspace{2mm}
\sigma(x,x) = \nu_1\,b; \\
\vspace{2mm}
\sigma(x,y) = \lambda\,b + \mu\,l;  \\
\vspace{2mm}
\sigma(y,y) = \nu_2\,b,
\end{array}\leqno{(2.1)}$$
where $\nu_1, \nu_2, \lambda, \mu$ are invariant functions, whose signs depend on the
pair $\{b, l\}$.

Hence the invariants $k$, $\varkappa$, and the Gauss curvature $K$ of $M^2$ are expressed
as follows:
$$k = - 4\nu_1\,\nu_2\,\mu^2, \quad \quad \varkappa = (\nu_1-\nu_2)\mu, \quad \quad
K = \nu_1\,\nu_2- (\lambda^2 + \mu^2).\leqno(2.2)$$

The normal mean curvature vector field $H$ of $M^2$ is
$H = \ds{\frac{\sigma(x,x) + \sigma(y,y)}{2} = \frac{\nu_1 + \nu_2}{2}\, b}$.
\vskip 2mm
Let $M^2$ be a surface parameterized by principal tangents and $g: X = \alpha z_u + \beta z_v$
be an arbitrary tangent of $M^2$. We call the function
$\ds{\nu_g = \frac{II(\alpha, \beta)}{I(\alpha, \beta)}}$ the \textit{normal curvature} of $g$.
Obviously, a tangent $g$ is asymptotic if and only if its normal curvature is zero.

The normal curvatures $\nu' = \ds{\frac{L}{E}}$ and $\nu'' = \ds{\frac{N}{G}}$ of
the principal tangents are said to be \textit{principal normal curvatures} of $M^2$.
If $g$ is an arbitrary tangent with normal curvature $\nu_g$, and $\varphi = \angle (g,z_u)$,
then the following Euler formula holds
$$\nu_g = \cos^2 \varphi \, \nu' + \sin^2 \varphi \, \nu''.$$
The invariants $k$ and $\varkappa$ of $M^2$ are expressed by the principal normal curvatures
$\nu'$ and $\nu''$ as follows:
$$k = \nu' \nu''; \qquad \varkappa = \frac{\nu' + \nu''}{2}. \leqno{(2.3)}$$
Hence, the invariants $k$ and $\varkappa$ of $M^2$ play the same role in the differential
geometry of surfaces in $\R^4$ as the Gaussian curvature and the mean curvature
in the classical differential geometry of surfaces in $\R^3$.
\vskip 1mm
As in the theory of surfaces in $\R^3$, we consider the indicatrix $\chi$ in the tangent space
$T_pM^2$ at an arbitrary point $p$ of $M^2$, defined by
$$\chi : \nu' X^2 + \nu'' Y^2 = \varepsilon, \qquad \varepsilon = \pm 1.$$

If $p$ is an elliptic point ($k > 0$), then the indicatrix $\chi$ is an ellipse.
The axes of $\chi$ are collinear with the principal directions at the point $p$,
and the lengths of the axes are
$\ds{\frac{2}{\sqrt{|\nu'|}}}$ and $\ds{\frac{2}{\sqrt{|\nu''|}}}$.

If $p$ is a  hyperbolic point ($k < 0$), then the indicatrix $\chi$ consists of two hyperbolas.
For the sake of simplicity we say that $\chi$ is a hyperbola. The axes of $\chi$ are collinear
with the principal directions, and the lengths of the axes are
$\ds{\frac{2}{\sqrt{|\nu'|}}}$ and $\ds{\frac{2}{\sqrt{|\nu''|}}}$.

If $p$ is a parabolic point ($k = 0$), then the indicatrix $\chi$
consists of two straight lines parallel to the principal direction with non-zero normal curvature.

The following statement holds good:

\begin{prop}\label{P:conjugacy}
Two tangents $g_1$ and $g_2$ are conjugate tangents of $M^2$ if and only if
$g_1$ and $g_2$ are conjugate with respect to the indicatrix $\chi$.
\end{prop}
\vskip 2mm
\section{Classes of surfaces characterized in terms of the tangent indicatrix and the normal curvature ellipse}
\vskip 2mm
Each surface $M^2$ in $\R^4$ satisfies the following inequality:
$$\varkappa^2 - k \geq 0.$$

The minimal surfaces in $\R^4$ are characterized by

\begin{prop}\label{P:minimal} \cite{GM1}
Let $M^2$ be a surface in $\R^4$ free of flat points. Then $M^2$ is minimal if and only if
$$\varkappa^2 - k = 0.$$
\end{prop}

The surfaces with flat normal connection are characterized by

\begin{prop}\label{P:Flat normal}
Let $M^2$ be a surface in $\R^4$ free of flat points. Then $M^2$ is a surface with
flat normal connection if and only if
$$\varkappa = 0.$$
\end{prop}
\vskip 2mm
We note that the condition $\varkappa = 0$ implies that $k<0$ and the surface $M^2$ has
two families of orthogonal asymptotic lines.
\vskip 2mm
Now we shall characterize the minimal surfaces and the surfaces with flat normal connection
in terms of the tangent indicatrix of the surface.

\begin{prop}\label{P:Minimal-circle}
Let $M^2$ be a surface in $\R^4$ free of flat points. Then $M^2$ is minimal if and only if
at each point of $M^2$ the tangent indicatrix $\chi$ is a circle.
\end{prop}

\noindent
{\it Proof:} Let $M^2$ be a surface in $\R^4$ free of flat points.
From equalities (2.3) it follows that
$$\varkappa^2 - k = \ds{\left(\frac{\nu' - \nu''}{2}\right)^2}.$$
Obviously $\varkappa^2 - k = 0$ if and only if $\nu' = \nu''$.  Applying
Proposition \ref{P:minimal}, we get that $M^2$ is minimal if and only if $\chi$ is a circle.
\qed

\begin{prop}\label{P:Flat normal-Lorentz circle}
Let $M^2$ be a surface in $\R^4$ free of flat points.  Then $M^2$ is a surface
of flat normal connection if and only if at each point of $M^2$ the tangent indicatrix
$\chi$ is a rectangular hyperbola (a Lorentz circle).
\end{prop}

\noindent
{\it Proof:} Let $M^2$ be a surface in $\R^4$ free of flat points. From (2.3)
it follows that $\varkappa = 0$ if and only if $\nu'' = -\nu'$.

If $M^2$ is a surface with flat normal connection, then $k < 0$, and hence
$\chi$ is a hyperbola. From $\nu'' = -\nu'$ it follows that the semi-axes of
$\chi$ are equal to $\ds{\frac{1}{\sqrt{|\nu'|}}}$, i.e. $\chi$ is a rectangular hyperbola.

Conversely, if $\chi$ is a rectangular hyperbola, then $\nu'' = -\nu'$,
which implies that $M^2$ is a surface with flat normal connection.
\qed
\vskip 2mm
The minimal surfaces and the surfaces with flat normal connection can also be characterized
in terms of the ellipse of normal curvature.

Let us recall that the \textit{ellipse of normal curvature} at a point $p$ of a surface
$M^2$ in $\R^4$ is the ellipse in the normal space at the point $p$ given by
$\{\sigma(x,x): \, x \in T_pM^2, \, g(x,x) = 1\}$ \cite{MW1, MW2}.
Let $\{x,y\}$ be an orthonormal base of the tangent space $T_pM^2$ at $p$. Then, for any
$v = \cos \psi \, x + \sin \psi \, y$, we have
$$\sigma(v, v) = H + \ds{\cos 2\psi  \, \frac{\sigma(x,x) - \sigma(y,y)}{2}
+ \sin 2 \psi  \, \sigma(x,y)},$$
where $H = \ds{\frac{\sigma(x,x) + \sigma(y,y)}{2}}$ \, is the mean curvature vector of $M^2$
at $p$. So, when $v$ goes once around the unit tangent circle,
the vector $\sigma(v,v)$ goes twice around the ellipse centered at $H$. The vectors
$\ds{\frac{\sigma(x,x) - \sigma(y,y)}{2}}$ \, and $\sigma(x,y)$ determine conjugate
directions of the ellipse.

A surface $M^2$ in $\R^4$ is called \textit{super-conformal} \cite{DT} if at any point of $M^2$
the ellipse of curvature is a circle. In \cite{DT} it is given an explicit construction of any
simply connected super-conformal surface in $\R^4$ that is free of minimal and flat points.
\vskip 2mm
Obviously, $M^2$ is minimal if and only if
for each point $p \in M^2$ the ellipse of curvature is centered at $p$.

The  minimal surfaces in $\R^4$ are divided into two subclasses:
\begin{itemize}
\item
the subclass of minimal super-conformal surfaces, characterized by the condition that
the ellipse of curvature is a circle;

\item subclass of minimal surfaces of general type, characterized by the condition that
the ellipse of curvature is not a circle.
\end{itemize}

In  \cite{GM2} it is proved that on any minimal surface $M^2$ the Gauss
curvature $K$ and the normal curvature $\varkappa$ satisfy the following inequality
$$K^2-\varkappa^2\geq 0.$$
The two subclasses of minimal surfaces are characterized in terms of the invariants $K$
and  $\varkappa$ as follows:
\begin{itemize}
\item the class of minimal super-conformal surfaces is characterized by $K^2 - \varkappa^2 =0$;
\item the class of minimal surfaces of general type is characterized by $K^2-\varkappa^2>0$.
\end{itemize}

The class of minimal super-conformal surfaces in $\R^4$ is locally equivalent to the class
of holomorphic curves in $\C^2 \equiv \R^4$.
\vskip 2mm
The surfaces with flat normal connection are characterized in terms of the ellipse of
normal curvature as follows

\begin{prop}\label{P:flat-normal-ellipse}
Let $M^2$ be a surface in $\R^4$ free of flat points. Then $M^2$ is a surface with
flat normal connection if and only if for each point $p \in M^2$ the ellipse of normal
curvature is a line segment, which is not collinear with the mean curvature vector field.
\end{prop}

\noindent
{\it Proof:}
In \cite{A} it is proved that the curvature of the normal connection $\varkappa$
of a surface $M^2$ in $\R^4$ is the Gauss torsion $\varkappa_G$ of $M^2$.
The notion of the Gauss torsion is introduced by \'E. Cartan \cite{C} for a $p$-dimensional
submanifold of an $n$-dimensional Riemannian manifold and is given by the Euler curvatures.
In case of a 2-dimensional surface $M^2$ in $\R^4$ the Gauss torsion at a point $p \in M^2$
is equal to $2 ab$, where $a$ and $b$ are the semi-axis of the ellipse of normal curvature at $p$.
Hence, $\varkappa = 0$ if and only if the ellipse of curvature is a line segment.

Let $M^2$ be a surface with flat normal connection, i.e. $\varkappa = 0$, $k \neq 0$.
From (2.2) it follows, that $\nu_1 = \nu_2$. Further, equalities (2.1) imply that for each
$v = \cos \psi \, x + \sin \psi \, y$, we have
$\sigma(v, v) = H + \sin 2 \psi (\lambda\,b + \mu\,l)$.
So, when $v$ goes once around the unit tangent circle,
the vector $\sigma(v,v)$ goes twice along the line segment collinear with
$\lambda\,b + \mu\,l$ and centered at $H$.
The mean curvature vector field is $H = \nu_1\, b$.  Since $k \neq 0$ then
$\mu \neq 0$, and the line segment is not collinear with $H$. \qed
\vskip 2mm
In case of $\lambda = 0$ the mean curvature vector field $H$ is orthogonal to
the line segment, while in case of $\lambda \neq 0$ the mean curvature vector field $H$
is not orthogonal to the line segment. The length $d$ of the line segment is
$$d = \sqrt{\lambda^2 + \mu^2} = \sqrt {H^2 - K}.$$
So, there arises a subclass of surfaces with flat normal connection, characterized by the
conditions:
$$ K = 0  \quad {\rm or} \quad d = \Vert H \Vert.$$

Proposition \ref{P:Flat normal-Lorentz circle} and Proposition \ref{P:flat-normal-ellipse}
give us the following

\begin{cor}\label{C:Lorentz circle-normal ellipse}
Let $M^2$ be a surface in $\R^4$ free of flat points.  Then the tangent indicatrix
$\chi$ is a rectangular hyperbola (a Lorentz circle) if and only if the ellipse of
normal curvature is a line segment, which is not collinear with the mean curvature
vector field.
\end{cor}

\section{Examples of surfaces with flat normal connection} \label{S:Examples}

In this section we construct a family of surfaces with flat normal connection lying
on a standard rotational hypersurface in $\R^4$ .

Let $\{e_1, e_2, e_3, e_4\}$ be the standard orthonormal frame in $\R^4$, and $S^2(1)$
be a 2-dimensional sphere in $\R^3 = \span \{e_1, e_2, e_3\}$, centered at the origin $O$.
We consider a smooth curve $c: l = l(v), \, v \in J, \,\, J \subset \R$  on $S^2(1)$,
parameterized by the arc-length ($l'^2(v) = 1$). We denote $t = l'$ and consider the
moving frame field $\span \{t(v), n(v), l(v)\}$ of the curve $c$ on $S^2(1)$.
With respect to this orthonormal frame field the following Frenet formulas hold good:
$$\begin{array}{l}
\vspace{2mm}
l' = t;\\
\vspace{2mm}
t' = \kappa \,n - l;\\
\vspace{2mm}
n' = - \kappa \,t,
\end{array} \leqno{(4.1)}$$
where $\kappa$ is the spherical curvature of $c$.

Let $f = f(u), \,\, g = g(u)$ be smooth functions, defined in an interval $I \subset \R$,
such that $\dot{f}^2(u) + \dot{g}^2(u) = 1, \,\, u \in I$.
Now we construct a surface $M^2$ in $\R^4$ in the following way:
$$M^2: z(u,v) = f(u) \, l(v) + g(u)\, e_4, \quad u \in I, \, v \in J. \leqno{(4.2)}$$

The surface $M^2$ lies on the rotational hypersurface $M^3$ in $\R^4$ obtained by the
rotation of the meridian curve $m: u \rightarrow (f(u), g(u))$ around the
$Oe_4$-axis in $\R^4$. Since $M^2$ consists of meridians of $M^3$, we call $M^2$
a \textit{meridian surface}.

The tangent space of $M^2$ is spanned by the vector fields:
$$\begin{array}{l}
\vspace{2mm}
z_u = \dot{f} \,l + \dot{g}\,e_4;\\
\vspace{2mm}
z_v = f\,t,
\end{array}$$
and hence the coefficients of the first fundamental form of $M^2$ are
$E = 1; \,\, F = 0; \,\, G = f^2(u)$.
Taking into account (4.1), we calculate the second partial derivatives of $z(u,v)$:
$$\begin{array}{l}
\vspace{2mm}
z_{uu} = \ddot{f} \,l + \ddot{g}\,e_4;\\
\vspace{2mm}
z_{uv} = \dot{f}\,t;\\
\vspace{2mm}
z_{vv} = f \kappa \,n - f\,l.
\end{array}$$
Let us denote $x = z_u,\,\, y = \ds{\frac{z_v}{f} = t}$ and consider the following
orthonormal normal frame field of $M^2$:
$$n_1 = n(v); \qquad n_2 = - \dot{g}(u)\,l(v) + \dot{f}(u) \, e_4.$$
Thus we obtain a positive orthonormal frame field $\{x,y, n_1, n_2\}$ of $M^2$.
If we denote by $\kappa_m$ the curvature of the meridian curve $m$, i.e.
$\kappa_m (u)= \dot{f}(u) \ddot{g}(u) - \dot{g}(u) \ddot{f}(u) =
\ds{\frac{- \ddot{f}(u)}{\sqrt{1 - \dot{f}^2(u)}}}$,
then we get the following derivative formulas of $M^2$:
$$\begin{array}{ll}
\vspace{2mm} \nabla'_xx = \qquad \qquad \qquad \qquad \kappa_m\,n_2;
& \qquad
\nabla'_x n_1 = 0;\\
\vspace{2mm}
\nabla'_xy = 0;  & \qquad
\nabla'_y n_1 = \ds{\quad \quad \quad - \frac{\kappa}{f}\,y};\\
\vspace{2mm}
\nabla'_yx = \quad\quad \quad\ds{\frac{\dot{f}}{f}}\,y;  & \qquad
\nabla'_x n_2 = - \kappa_m \,x;\\
\vspace{2mm}
\nabla'_yy = \ds{- \frac{\dot{f}}{f}\,x \quad\quad + \frac{\kappa}{f}\,n_1 + \frac{\dot{g}}{f} \, n_2};
& \qquad
\nabla'_y n_2 = \ds{ \quad \quad \quad - \frac{\dot{g}}{f}\,y}.
\end{array} \leqno{(4.3)}$$

The coefficients of the second fundamental form of $M^2$ are
$L = N = 0, \,\, M = - \kappa_m(u) \, \kappa(v)$.
Taking into account (4.3), we find the invariants $k$, $\varkappa$, $K$:
$$k = - \frac{\kappa_m^2(u) \, \kappa^2(v)}{f^2(u)}; \qquad \varkappa = 0;
\qquad K = \frac{\kappa_m (u)\, \dot{g}(u)}{f(u)}.  \leqno{(4.4)}$$

The equality $\varkappa = 0$ implies that $M^2$ is a surface with flat normal connection.

The mean curvature vector field $H$ is given by
$$H = \frac{\kappa}{2f}\, n_1 + \frac{\dot{g} + f \kappa_m}{2f} \, n_2. \leqno{(4.5)}$$

There are three main classes of meridian surfaces:
\vskip 2mm
I. $\kappa = 0$, i.e. the curve $c$ is a great circle on $S^2(1)$.
In this case $n_1 = const$, and $M^2$ is a planar surface
lying in the constant 3-dimensional space spanned by $\{x, y, n_2\}$.
Particularly, if in addition $\kappa_m = 0$, i.e. the meridian curve lies on a straight line,
then $M^2$ is a developable surface in the 3-dimensional space $\span \{x, y, n_2\}$.

\vskip 2mm
II. $\kappa_m = 0$, i.e. the meridian curve is part of a straight line. In such case
$k = \varkappa = K = 0$, and $M^2$ is a developable ruled surface. If in addition
$\kappa = const$, i.e. $c$ is a circle on $S^2(1)$, then
$M^2$ is a developable ruled surface in a 3-dimensional space.
If $\kappa \neq const$, i.e. $c$ is not a circle on $S^2(1)$, then
$M^2$ is a developable ruled surface in $\R^4$.

\vskip 2mm
III. $\kappa_m \, \kappa \neq 0$, i.e. $c$ is not a great circle on $S^2(1)$, and
$m$ is not a straight line. In this general case the invariant function $k<0$,
which implies that there exist two systems of asymptotic lines on $M^2$. The
parametric lines of $M^2$ given by (4.2) are orthogonal and asymptotic.

\vskip 2mm
Let $M^2$ be a meridian surface of the general class. Now we are going to find
the meridian surfaces with:
\begin{itemize}
\item
constant Gauss curvature $K$;
\item
constant mean curvature;
\item
constant invariant function $k$.
\end{itemize}

\begin{prop}
Let $M^2$ be a meridian surface in $\R^4$. Then $M^2$ has constant non-zero Gauss curvature $K$
if and only if the meridian $m$ is given by
$$\begin{array}{ll}
\vspace{2mm}
f(u) = \alpha \cos \sqrt{K} u + \beta \sin \sqrt{K} u, & \quad K >0;\\
\vspace{2mm}
f(u) = \alpha \cosh \sqrt{-K} u + \beta \sinh \sqrt{-K} u, & \quad K <0,
\end{array}$$
where $\alpha$ and $\beta$ are constants.
\end{prop}

\noindent
{\it Proof:}
Using (4.4) and $\dot f^2+\dot g^2= 1$, we obtain that $M^2$ has constant Gauss curvature
$K \neq 0$ if and only if the meridian $m$ satisfies the following differential equation
$$\ddot{f}(u) + K f(u) = 0.$$
The general solution of the above equation is given by
$$\begin{array}{ll}
\vspace{2mm}
f(u) = \alpha \cos \sqrt{K} u + \beta \sin \sqrt{K} u, & \quad \textrm{in case} \quad K >0;\\
\vspace{2mm}
f(u) = \alpha \cosh \sqrt{-K} u + \beta \sinh \sqrt{-K} u, & \quad \textrm{in case}  \quad K <0,
\end{array}$$
where $\alpha$ and $\beta$ are constants.
The function $g(u)$ is determined by $\dot{g}(u) = \sqrt{1 - \dot{f}^2(u)}$.

\qed

\vskip 3mm
The equality (4.5) implies that the mean curvature of $M^2$ is given by
$$|| H || = \ds{\sqrt{\frac{\kappa^2(v) + \left(\dot{g}(u)
+ f(u) \kappa_m(u)\right)^2}{4f^2(u)}}}. \leqno{(4.6)}$$

The  meridian surfaces with constant mean curvature (CMC meridian surfaces) are described in

\begin{prop}
Let $M^2$ be a meridian surface in $\R^4$. Then $M^2$ has constant mean curvature
$|| H || = a = const$, $a \neq 0$
if and only if the curve $c$  on $S^2(1)$ is a circle with constant spherical curvature
$\kappa = const = b, \; b \neq 0$,
and the meridian $m$ is determined by the following differential equation:
$$\left( 1 - \dot{f}^2 - f \ddot{f}\right)^2 = (1 - \dot{f}^2) (4 a^2 f^2 - b^2).$$
\end{prop}

\noindent
{\it Proof:}
From (4.6) it follows that $||H|| = a$ if and only if
$$\kappa^2(v) = 4 a^2 f^2(u) - (\dot{g}(u) + f(u)  \kappa_m(u))^2,$$
which implies
$$\begin{array}{l}
\vspace{2mm}
\kappa = const = b, \; b \neq 0;\\
\vspace{2mm}
4 a^2 f^2(u) - (\dot{g}(u) + f(u)  \kappa_m(u))^2 = b^2.
\end{array} \leqno{(4.7)}$$
The first equality of (4.7) implies that the spherical curve $c$ has constant
spherical curvature $\kappa = b$,  i.e. $c$ is a circle.
Using that $\dot{f}^2 + \dot{g}^2 = 1$, and $\kappa_m = \dot{f} \ddot{g} - \dot{g} \ddot{f}$
we calculate that
$\dot{g} + f \kappa_m = \ds{\frac{1- \dot{f}^2 - f \ddot{f}}{\sqrt{1 - \dot{f}^2}}}.$
Hence, the second equality of (4.7) gives the following
differential equation for the meridian $m$:
$$\left( 1 - \dot{f}^2 - f \ddot{f}\right)^2 = (1 - \dot{f}^2) (4 a^2 f^2 - b^2). \leqno{(4.8)}$$

Further, if we set $\dot{f} = y(f)$ in equation (4.8), we obtain that the function
$y = y(t)$ is a solution of the following differential equation
$$1 - y^2 - \frac{t}{2}(y^2)' = \sqrt{1 - y^2} \sqrt{4 a^2 t^2 - b^2}.$$
The general solution of the above equation is given by
$$y(t) = \sqrt{1 - \frac{1}{t^2}\left( C + \frac{t}{2} \sqrt{4 a^2 t^2 - b^2} -
\frac{b^2}{4a} \ln |2at + \sqrt{4 a^2 t^2 - b^2}| \right)^2}; \qquad C = const. \leqno{(4.9)}$$
The function $f(u)$ is determined by $\dot{f} = y(f)$ and (4.9).
The function $g(u)$ is defined by $\dot{g}(u) = \sqrt{1 - \dot{f}^2(u)}$.
\qed
\vskip 2mm
At the end of this section we shall find the meridian surfaces with constant invariant $k$.

\begin{prop}
Let $M^2$ be a meridian surface in $\R^4$. Then $M^2$ has a constant invariant
$k = const = - a^2, \; a\neq 0$ if and only if the curve $c$  on $S^2(1)$ is
a circle with spherical curvature $\kappa = const = b, \; b \neq 0$,
and the meridian $m$ is determined by the following differential equation:
$$\ddot{f}(u) = \mp \frac{a}{b}\,f(u) \sqrt{1 - \dot{f}^2(u)}.$$
\end{prop}

\noindent
{\it Proof:}
Using (4.4) we obtain that $k = const = - a^2, \; a \neq 0$
if and only if $\kappa^2(v) \kappa_m^2(u) = a^2 f^2(u)$. Hence,
$$\kappa(v) = \pm \, a \, \frac{f(u)}{\kappa_m(u)}.$$
The last equality implies
$$\begin{array}{l}
\vspace{2mm}
\kappa = const = b, \; b \neq 0;\\
\vspace{2mm}
 \pm \, a \, \ds{\frac{f(u)}{\kappa_m(u)}} = b.
\end{array} \leqno{(4.10)}$$
The first equality of (4.10) implies that the spherical curve $c$ has constant spherical
curvature $\kappa = b$,  i.e. $c$ is a circle. The second equality of (4.10) gives
the following differential equation for the function $f(u)$:
$$\frac{\ddot{f}(u)}{\sqrt{1 - \dot{f}^2(u)}} = \mp \frac{a}{b}\,f(u).  \leqno{(4.11)}$$
Again setting $\dot{f} = y(f)$ in equation (4.11), we obtain that the function
$y = y(t)$ is a solution of the following differential equation

$$\frac{y y'}{\sqrt{1 - y^2}} = \mp \frac{a}{b}\,t.$$
The general solution of the above equation is given by
$$y(t) = \sqrt{1 - \left( C \pm \frac{a}{b}\, \frac{t^2}{2} \right)^2}; \qquad C = const.  \leqno{(4.12)}$$
The function $f(u)$ is determined by $\dot{f} = y(f)$ and (4.12).
The function $g(u)$ is defined by $\dot{g}(u) = \sqrt{1 - \dot{f}^2(u)}$.
\qed
\vskip 2mm
\section{Examples of surfaces consisting of parabolic points} \label{S:Parabolic points}
\vskip 2mm
In this section we shall find the generalized (in the sense of C. Moore) rotational surfaces in $\R^4$,
consisting of parabolic points.
\vskip 1mm
We consider a surface $M^2$ in $\R^4$ given by
$$z(u,v) = \left( f(u) \cos\alpha v, f(u) \sin \alpha v, g(u) \cos \beta v, g(u) \sin \beta v \right);
\quad u \in J \subset \R, \,\,  v \in [0; 2\pi), \leqno{(5.1)}$$
where $f(u)$ and $g(u)$ are smooth functions, satisfying
$\alpha^2 f^2(u)+ \beta^2 g^2(u) > 0 , \,\, f'\,^2(u)+ g'\,^2(u) > 0,\, u \in J$,
and $\alpha, \beta$ are positive constants.

Each parametric curve $u = u_0 = const$ of $M^2$ is given by
$$c_v: z(v) = \left( a \cos \alpha v, a \sin \alpha v, b \cos \beta v, b \sin \beta v \right);
\quad a = f(u_0), \,\, b = g(u_0)$$ and
its Frenet curvatures are
$$\varkappa_{c_v} = \ds{\sqrt{\frac{a^2 \alpha^4 + b^2 \beta^4}{a^2 \alpha^2 + b^2 \beta^2}}}; \quad
\tau_{c_v} = \ds{\frac{ab \alpha \beta (\alpha^2 - \beta^2)}
{\sqrt{a^2 \alpha^4 + b^2 \beta^4}\sqrt{a^2 \alpha^2 + b^2 \beta^2}}}; \quad
\sigma_{c_v} = \ds{\frac{\alpha \beta \sqrt{a^2 \alpha^2 + b^2 \beta^2}}
{\sqrt{a^2 \alpha^4 + b^2 \beta^4}}}.$$
Hence, in case of $\alpha \neq \beta$ each parametric curve $u = const$ is a curve in $\R^4$
with constant curvatures, and in case of $\alpha = \beta$ each parametric curve
$u = const$ is a circle.

Each parametric curve $v = v_0 = const$ of $M^2$ is given by
$$c_u: z(u) = \left(\, A_1 f(u), A_2 f(u), B_1 g(u), B_2 g(u) \, \right),$$
where  $A_1 = \cos \alpha v_0, \, A_2 = \sin \alpha v_0, \,B_1 = \cos \beta v_0,
\,B_2 = \sin \beta v_0$.
The Frenet curvatures of $c_u$ are expressed as follows:
$$\varkappa_{c_u} = \ds{\frac{|g' f'' - f' g''|}{(\sqrt{f'\,^2 + g'\,^2})^3}}; \quad \tau_{c_u} = 0.$$
Hence, $c_u$ is a plane curve with curvature $\varkappa_{c_u} = \ds{\frac{|g' f'' - f' g''|}
{(\sqrt{f'\,^2 + g'\,^2})^3}}$.
So, for each $v = const$ the parametric curves $c_u$ are congruent in $\R^4$.
We call these curves \textit{meridians} of $M^2$.

Considering general rotations in $\R^4$, C. Moore introduced general rotational surfaces
\cite{M} \, (see also \cite{MW1, MW2}). The surface $M^2$, given by (5.1) is a general rotational surface
whose meridians lie in two-dimensional planes.

The tangent space of $M^2$ is spanned by the vector fields
$$\begin{array}{l}
\vspace{2mm}
z_u = \left(f' \cos \alpha v, f' \sin \alpha v, g' \cos \beta v, g' \sin \beta v \right);\\
\vspace{2mm}
z_v = \left( - \alpha f \sin \alpha v, \alpha f \cos \alpha v, -
\beta g \sin \beta v, \beta g \cos \beta v \right).
\end{array}$$
Hence, the coefficients of the first fundamental form are
$E = f'\,^2(u)+ g'\,^2(u); \,\, F = 0; \,\, G = \alpha^2 f^2(u)+ \beta^2 g^2(u)$
and $W = \sqrt{(f'\,^2 + g'\,^2)(\alpha^2 f^2 + \beta^2 g^2)}$. We consider the
following orthonormal tangent frame field
$$\begin{array}{l}
\vspace{2mm}
x = \ds{\frac{1}{\sqrt{f'\,^2 + g'\,^2}}\left(f' \cos \alpha v,
f' \sin \alpha v, g' \cos \beta v, g' \sin \beta v \right)};\\
\vspace{2mm}
y = \ds{\frac{1}{\sqrt{\alpha^2 f^2 + \beta^2 g^2}}\left( - \alpha f \sin \alpha v,
\alpha f \cos \alpha v, - \beta g \sin \beta v, \beta g \cos \beta v \right)}.
\end{array}$$
The second partial derivatives of $z(u,v)$ are expressed as
follows
$$\begin{array}{l}
\vspace{2mm}
z_{uu} = \left(f'' \cos \alpha v, f'' \sin \alpha v, g'' \cos \beta v, g'' \sin \beta v \right);\\
\vspace{2mm}
z_{uv} = \left(- \alpha f' \sin \alpha v, \alpha f' \cos \alpha v,
- \beta g' \sin \beta v, \beta g' \cos \beta v \right);\\
\vspace{2mm} z_{vv} = \left(- \alpha ^2 f \cos \alpha v,
- \alpha ^2 f \sin \alpha v, - \beta ^2 g \cos \beta v, -
\beta ^2 g \sin \beta v \right).
\end{array}$$

Now let us consider the following orthonormal normal frame field
$$\begin{array}{l}
\vspace{2mm}
n_1 = \ds{\frac{1}{\sqrt{f'\,^2 + g'\,^2}}\left(g' \cos \alpha v,
g' \sin \alpha v, - f' \cos \beta v, - f' \sin \beta v \right)};\\
\vspace{2mm} n_2 = \ds{\frac{1}{\sqrt{\alpha^2 f^2 + \beta^2 g^2}}\left( - \beta g \sin \alpha
v, \beta g \cos \alpha v, \alpha f \sin \beta v, - \alpha f \cos \beta v \right)}.
\end{array}$$
It is easy to verify that $\{x, y, n_1, n_2\}$ is a positive
oriented orthonormal frame field in $\R^4$.

We calculate the functions $c_{ij}^k, \,\, i,j,k = 1,2$:
$$\begin{array}{ll}
\vspace{2mm} c_{11}^1 = g(z_{uu}, n_1) = \ds{\frac{g'
f'' - f' g''}{\sqrt{f'\,^2 + g'\,^2}}};
\quad \quad & c_{11}^2 = g(z_{uu}, n_2) = 0;\\
\vspace{2mm} c_{12}^1 = g(z_{uv}, n_1) = 0;
\quad \quad & c_{12}^2 = g(z_{uv}, n_2) = \
ds{\frac{\alpha \beta (g f' - f g')}{\sqrt{\alpha^2 f^2 + \beta^2 g^2}}};\\
\vspace{2mm} c_{22}^1 = g(z_{vv}, n_1) =
\ds{\frac{\beta^2 g f' - \alpha^2 f g'}{\sqrt{f'\,^2 + g'\,^2}}};
\quad \quad & c_{22}^2 = g(z_{vv}, n_2) = 0.
\end{array} $$
Therefore the coefficients $L$, $M$ and $N$ of the second
fundamental form of $M^2$ are expressed as follows:
$$L = \ds{\frac{2 \alpha \beta (g f' - f g') (g' f'' - f' g'')}
{(\alpha^2 f^2 + \beta^2 g^2) (f'\,^2 + g'\,^2)}}; \qquad M = 0; \qquad
N = \ds{\frac{- 2\alpha \beta (g f' - f g') (\beta^2 g f' - \alpha^2 f g')}
{(\alpha^2 f^2 + \beta^2 g^2) (f'\,^2 + g'\,^2)}}.$$

Consequently, the invariants $k$, $\varkappa$ and $K$ of $M^2$
are:

$$k = \ds{\frac{- 4 \alpha^2 \beta^2 (g f' - f g')^2 (g' f'' - f' g'')
(\beta^2 g f' - \alpha^2 f g')}{(\alpha^2 f^2 + \beta^2 g^2)^3 (f'\,^2 + g'\,^2)^3}};$$

$$\varkappa =  \ds{\frac{\alpha \beta (g f' - f g')}
{(\alpha^2 f^2 + \beta^2 g^2)^2 (f'\,^2 + g'\,^2)^2} \,
\left( (\alpha^2 f^2 + \beta^2  g^2)(g' f'' - f' g'') -
(f'\,^2 + g'\,^2) (\beta^2 g f' - \alpha^2 f g')
\right)};$$

$$K =  \ds{\frac{(\alpha^2 f^2 + \beta^2  g^2)(\beta^2 g f' - \alpha^2 f g')
(g' f'' - f' g'') - \alpha^2 \beta^2 (f'\,^2 + g'\,^2) (g f' - f g')^2}
{(\alpha^2 f^2 + \beta^2 g^2)^2 (f'\,^2 + g'\,^2)^2} \,}.$$
\vskip 2mm
Now we shall find the generalized rotational surfaces with $k = 0$.
Without loss of generality we assume that
the meridian $m$ is defined by $f = u; \,\, g = g(u)$. Then
$$k = \ds{\frac{ 4 \alpha^2 \beta^2 (g  - u g')^2 g'' (\beta^2 g  - \alpha^2 u g')}
{(\alpha^2 u^2 + \beta^2 g^2)^3 (1 + g'\,^2)^3}};$$

The invariant $k$ is zero in the following three cases:
\vskip 2mm
1. $g(u) = a\,u,\; a = const \neq 0$. In that case $k = \varkappa = K = 0$,
and $M^2$ is a developable surface in $\R^4$.
\vskip 2mm
2. $g(u) = a\,u + b ,\; a= const \neq 0, b = const \neq 0$. In this case $k = 0$,
but $\varkappa \neq 0$, $K \neq 0$.
Consequently, $M^2$ is a non-developable ruled surface in $\R^4$.

\vskip 2mm
3. $g(u) = \ds{ c\,u^{\frac{\beta^2}{\alpha^2}}}, \; c = const \neq 0$.
In case of $\alpha \neq \beta$ we get $k = 0$,
and the invariants $\varkappa$ and $K$ are given by
$$\varkappa = \ds{\frac{c^2 \beta^3 (\beta^2 - \alpha^2 )^2 u^{2 \frac{\beta^2-\alpha^2}{\alpha^2}}}
{\alpha^5 \left(\alpha^2 u^2 + \beta^2 c^2 u^{2 \frac{\beta^2}{\alpha^2}}\right)
\left(1 + c^2 \frac{\beta^4}{\alpha^4} u^{2\frac{\beta^2-\alpha^2}{\alpha^2} }\right)^2}};$$
$$K = \ds{- \frac{c^2 \beta^2 (\beta^2 - \alpha^2 )^2 u^{2 \frac{\beta^2}{\alpha^2}}}
{\alpha^2 \left(\alpha^2 u^2 + \beta^2 c^2 u^{2 \frac{\beta^2}{\alpha^2}}\right)^2
\left(1 + c^2 \frac{\beta^4}{\alpha^4} u^{2\frac{\beta^2-\alpha^2}{\alpha^2}}\right)}}.$$
Hence, $\varkappa \neq 0$, $K \neq 0$.
In this case the parametric lines $u =const$ and $v = const$ are not straight lines.
This is a non-trivial example of generalized rotational surfaces with $k = 0$.
\vskip 3mm
\textbf{Acknowledgements:} The second author is partially supported by "L. Karavelov" Civil
Engineering Higher School, Sofia, Bulgaria under Contract No 10/2009.

\end{document}